\documentclass [12pt]{article}
\usepackage{graphicx,amssymb,amsfonts,latexsym,amsmath,amsthm,times}
\usepackage{epsfig}
\usepackage{fancyhdr}
\usepackage{color}
\usepackage[english]{babel}
\usepackage{dsfont} \usepackage{color}
\usepackage{graphicx} \usepackage{float}
\usepackage{multirow}
\newcommand{\edp}[3]{\begin{equation}\label{#1}\left\lbrace\begin{array}{#2}
#3
\end{array}\right.\end{equation}}

\numberwithin{equation}{section}

\begin{document}




\thispagestyle{plain}

\vspace*{2cm} \normalsize \begin{center}
{\Large \bf A mathematical model of systemic inhibition of angiogenesis in metastatic development}
\end{center}
\vspace*{1cm}

\centerline{\bf S. Benzekry$^{1}$\footnote{Corresponding author. E-mail: sebastien.benzekry@inria.fr},  A. Gandolfi$^{2}$, P. Hahnfeldt$^3$}

\vspace*{0.5cm}

\centerline{$^1$ Inria team MC2, Institut de Math\'{e}matiques de Bordeaux, Bordeaux, France}
\centerline{$^2$ Istituto di Analisi dei Sistemi ed Informatica ``Antonio Ruberti'' - CNR, Roma, Italy}
\centerline{$^3$ Center of Cancer Systems Biology, GRI, Tufts University School of Medicine, Boston, 02142, USA}


\vspace*{1cm}

\noindent {\bf Abstract.}
We present a mathematical model describing the time development of a population of tumors subject to mutual angiogenic inhibitory signaling. Based on biophysical derivations, it describes organism-scale population dynamics under the influence of three processes: birth (dissemination of secondary tumors), growth and inhibition (through angiogenesis). The resulting model is a nonlinear partial differential transport equation with nonlocal boundary condition. The nonlinearity stands in the velocity through a nonlocal quantity of the model (the total metastatic volume). The asymptotic behavior of the model is numerically investigated and reveals interesting dynamics ranging from convergence to a steady state to bounded non-periodic or periodic behaviors, possibly with complex repeated patterns.  Numerical simulations are performed with the intent to theoretically study the relative impact of potentiation or impairment of each process of the birth/growth/inhibition balance. Biological insights on possible implications for the phenomenon of ``cancer without disease" are also discussed.

\vspace*{0.5cm}

\noindent {\bf Key words:} Nonlinear renewal equation; Cancer modeling; Metastasis; Angiogenesis inhibition;

\vspace*{0.5cm}
\noindent {\bf AMS subject classification:} 35B40, 35Q92, 92C50


\vspace*{1cm}

\newpage
\tableofcontents

\setcounter{equation}{0}


\section{Introduction}
Metastasis, the process by which secondary tumors are shed from a primary lesion to colonize local or distant sites, is a complex process that is responsible for a very large proportion (90\%) of deaths by a cancer disease \cite{Fidler2003}. Despite significant efforts in the last decade to strengthen our understanding of metastatic cancer biology, global mechanisms are still poorly understood \cite{Hanahan2011}. 

Among the wide variety of topics to be addressed, we focus here on signaling interactions within a population of tumors that impact on the global dynamics of the system, at the organism scale. Indeed, it has been long known that several tumors simultaneously growing within the same host influence each other, mostly in an inhibitory fashion. Most of the experimental studies were conducted in two-tumors experimental systems where impaired growth was observed for a second inoculum in the presence of a pre-existing tumor, a phenomenon that has been termed concomitant resistance \cite{Prehn1993} (see \cite{Chiarella2012} for a review). This distant inhibition not only occurs between two artificially implanted tumors but also between a primary tumor and its metastases \cite{Dewys1972a, Gorelik1978}. Such an observation is of fundamental importance for cancer biology as it impacts on the temporal development of the disease, but also has clinical implications in terms of metastatic dormancy \cite{Aguirre-Ghiso2007} and surgery \cite{Retsky2010}.

Several hypotheses have been proposed to explain this phenomenon, among which athrepsia (deprivation of nutriments for the second tumor due to monopolization of nutrients by the first one) or immune enhancement. Indeed, it is known \cite{Hanahan2011} that some tumors are immunogenic, i.e. that they provoke a hostile reaction from the immune system. This reaction could be triggered by the presence of a first tumor and suppress the growth of a subsequent implant. However, the importance of this last hypothesis was tempered in the 1980's when concomitant resistance was shown to occur in immunodeficient mice \cite{Gorelik1983}.

In the 1990's, Judah Folkman and his team put forward a novel hypothesis, based on their discovery of endogenous inhibitors of tumor neo-angiogenesis \cite{O'Reilly1994,O'Reilly1997}. Angiogenesis is the process of creation of new blood vessels from an existing vasculature and was shown to be of fundamental importance for a tumor's development \cite{Folkman1971} in order to ensure access to nutrients. Indeed, in the absence of angiogenesis, a neoplastic lesion cannot grow beyond 2-3 mm in diameter \cite{Folkman1971}. Observation that a tumor not only stimulates this process (through production of growth factors) but also regulates it by producing angiogenesis inhibitors opened the way to a new explanation of concomitant resistance. A primary tumor could distantly inhibit a secondary tumor through inhibitory angiogenic signaling, an hypothesis that is strengthened by the observation that angiogenic inhibitors have significant longer half-lives than angiogenic stimulators (of the order of hours for inhibitors such as angiostatin \cite{O'Reilly1994} or thrombospondin-1 \cite{Rofstad2001} and of the order of minutes for stimulating agents such as VEGF \cite{Folkman1995}). The formers are thus transferred to the central circulation and from there systemically distributed to distant sites. Considering the large and unequivocal body of support for the role angiogenesis inhibition plays in the maintenance of tumor dormancy \cite{Holmgren1995, O'Reilly1994, O'Reilly1997, Rofstad2001, Volpert1998, Sckell1998} we will focus on systemic inhibition of angiogenesis (SIA) as the major process for tumor-tumor interactions at the organism scale.

We present here a mathematical model for these interactions with the aim of qualitatively study the resulting nonlinear dynamics at the organism scale. The approach we use is based on previous work \cite{BenzekryJEE,BenzekryM2AN} on modeling of metastatic development, in the framework of structured population dynamics. 

\section{Model}
The main idea, originally introduced by Iwata et al. \cite{iwata}, and further studied in \cite{BBHV, Devys2009} is to represent the population of metastatic colonies as a density $\rho$ structured in tumor macroscopic traits (such as the volume in \cite{iwata}) and to derive a transport partial differential equation on $\rho$ that reflects mass conservation during the growth, endowed with a nonlocal boundary condition for dissemination (birth) of new metastases. 

A major limitation of this initial model for our purpose is that it does not take angiogenesis into account, although this is a fundamental process of tumor development. By combining the approach of \cite{iwata} with the model of \cite{Hahnfeldt1999} for tumor growth under angiogenic control, we developed in previous work a new model taking into account the angiogenic process in the growth of each tumor \cite{BenzekryM2AN,BenzekryJEE,BenzekryMMNP}. Combined to the fact that it is written at the level of the organism makes it an adapted framework for modeling SIA viewed as inhibiting interactions across a population of tumors. 

The model we propose here adds an inhibitory component to the two main processes of the previous model (angiogenic growth and metastatic dissemination). We add a new variable representing the circulating concentration of an endogenous angiogenesis inhibitor standing for all possible inhibitory molecules (examples being endostatin, angiostatin and thrombospondin-1).  It impacts on the growth of each tumor. As a general modeling principle, we want to be parsimonious and describe the major dynamics of the system with as few parameters as possible.

Tumors are seen as individuals whose state is described by two traits: volume $V$ and carrying capacity $K$, the latter representing environmental limitations that constrain growth of the tumor. This carrying capacity is here assumed to be equivalent to the vascular support. The primary tumor state is denoted $(V_p (t),K_p (t))$ while the model’s main variable for representation of the population of secondary tumors is denoted by $\rho(t,V,K)$ and stands for the physiologically structured density of metastases with volume $V$ and carrying capacity $K$ at time $t$. We assume that growth of the tumors is governed by a growth rate denoted $G(V,K;V_p(t),\rho)$, that spread of new metastases is driven by an volume-dependent dissemination rate denoted $\beta(V)$ and that the repartition of metastases at birth is given by a measure $N$. The precise expressions of these coefficients shall be described below. We consider some fixed final time $T$ and a physiological domain $\Omega=]V_0,+\infty[\times ]0,+\infty[$ where the distribution of metastases has its support, which means that only metastases having volume larger than some minimal volume $V_0$ and non-negative carrying capacity are considered biologically relevant. In the formula below, $\nu(\sigma)$ stands for the external normal to the boundary of the domain $\partial\Omega$. The notation $\partial\Omega^+$ stands for the subset of the boundary where the flux is pointing inward, i.e. where $G(\sigma;V_p(t),\rho(t))\cdot\nu(\sigma)<0$.  The function $\rho^0$ denotes the initial distribution of the metastatic colonies. Overall, the model writes
\edp{EqModel}{ll}{
\partial_t \rho (t,V,K) + {\rm div}\left(\rho(t,V,K)G(V,K;V_p(t),\rho(t))\right) = 0 & ]0,T[\times \Omega \\
-G(\sigma;V_p(t),\rho(t))\cdot \nu(\sigma)=N(\sigma)\left\lbrace\int_\Omega \beta(V)\rho(t,V,K)dVdK + \beta(V_p(t))\right\rbrace & ]0,T[\times \partial\Omega^+ \\
\rho(0,V,K)=\rho^0(V,K) & \Omega. 
}

\subsection{Tumor growth and systemic angiogenesis inhibition}
We assume that all the tumors (primary and secondaries) share the same growth model and parameters. We do so in order to reduce the number of parameters and focus on the dynamical properties of the system. Our population of tumors should thus be viewed as a population of identical entities in mutual interactions parallel to global development. The growth velocity of each tumor is given by a vector field $G(V,K;V_p(t),\rho(t))$. Following the approach of \cite{Hahnfeldt1999} we assume
$$G(V,K;V_p,\rho)=\left(
\begin{array}{c}
aV\ln\left(\frac VK\right) \\
Stim(V,K) - Inhib\left(V,K;V_p,\rho \right).
\end{array}
\right)$$
In the previous expression, the first line is the rate of change of the tumor volume $V$ (where $a$ is a constant parameter driving the proliferation kinetics) and the second line is the rate of change of the carrying capacity $K$. The main idea of this tumor growth model is to start from a gompertzian growth of the tumor volume (that could be replaced by any carrying capacity-like growth model, see  \cite{dOnofrio_Gandolfi}) and to assume that $K$ is a dynamical variable representing the tumor environment limitations (here limited to the vascular support) changing over time. The balance between a stimulation term $Stim(V,K)$ and an inhibition term denoted $Inhib(V,K;V_p(t),\rho(t,))$ (assumed here to depend on the global state of the system represented by the density $\rho$) governs the dynamics of $K$. For the stimulation term we follow \cite{Hahnfeldt1999} and assume 
$$Stim(V,K)=bV,$$
where the parameter $b$ is related to the concentration of angiogenic stimulating factors such as VEGF or bFGF.  This last quantity was derived to be constant in \cite{Hahnfeldt1999} from the consideration of very fast clearance of angiogenic stimulators \cite{Folkman1995}. This fact also strengthens our assumption of a local stimulatory term, as stimulating agents don't spread through the organism.

For the inhibition term, Hahnfeldt et al. \cite{Hahnfeldt1999} only considered a local inhibition coming from the tumor itself. Our main modeling novelty is to consider in addition a global inhibition coming from the release in the circulation of angiogenic inhibitors by the total (primary + secondary) population of tumors.  The following is an extension of the biophysical analysis performed in \cite{Hahnfeldt1999}. Let us consider a spherical tumor of radius R inside the host body. The host is represented, for simplicity, by a single compartment of volume $V_d$ in which concentrations are assumed spatially uniform. Let $n(r)$ be the inhibitor concentration inside the tumor at radial distance $r$. Let the intra-tumor clearance be zero \cite{Hahnfeldt1999}. At quasi steady state, $n(r)$ solves the following diffusion equation:
$$n''(r)+\frac{2n'(r)}{r}+\frac{p}{D^2} =0,$$
where $p$ is the inhibitor production rate and $D^2$ is the inhibitor diffusion coefficient. This equation is endowed with the boundary condition $n(R)=i(t;V_p(t),\rho(t))$, $i(t;V_p,\rho)$ being the systemic concentration of the inhibitor resulting from a primary tumor volume $V_p$ and secondary tumors density $\rho$ at time $t$. Solving this equation (using that $n(0)<+\infty$) we obtain
$$n(r)=i+\frac{p}{6D^2} (R^2-r^2 ).$$
From this expression we compute the mean inhibitor concentration in the tumor to obtain 
$$Inhib(V,K;V_p,ρ)=\hat{e}\left(i+\frac{p}{15D^2} R^2 \right)K=\hat{e}\left(i+\frac{p}{15D^2} \left(\frac{3}{4\pi}\right)^{2/3} V^{2/3}\right)K,$$
where $\hat{e}$ is a sensitivity coefficient. For $i(t;V_p(t),\rho(t))$, considering that the total flux of inhibitors produced by a tumor with volume $V$ is $pV$ and assuming that the inhibitor production rate is the same in all the tumors, we have
$$V_d  \frac{di}{dt}= pV_p(t)+\int_\Omega pV\rho(t,V,K)dVdK-kV_d i,$$
where $k$ is an elimination constant from the blood circulation. Defining $I(t;V_p(t),\rho(t))=V_d i(t;V_p(t),\rho(t))$, we get
\begin{equation}\label{EqInhib}
\frac{dI}{dt}= pV_p(t)+\int_\Omega pV\rho(t,V,K)dVdK-kI,
\end{equation}
endowed with zero initial condition ($I(0)=0$).

Overall, the explicit expression of the metastases growth rate is given by
\begin{equation}\label{EqGrowthRate}
G(V,K;V_p,\rho)=\left(
\begin{array}{c}
aV\ln\left(\frac VK\right) \\
bV-dV^{2/3}K-eI(V_p,\rho)K
\end{array}
\right),
\end{equation}
where $e:=\frac{\hat{e}}{V_d}$  and 
\begin{equation}\label{Eqd}
d:= e V_d  \frac{p}{15D^2} \left(\frac{3}{4\pi}\right)^{2/3}.
\end{equation} 
Note that we retrieve here the local term $dV^{2/3}$ from the analysis of \cite{Hahnfeldt1999}. Our analysis results in addition of a global term $eI$ for the effect of systemic inhibition of angiogenesis. 

For the primary tumor, we assume the growth velocity, hence the dynamics of $(V_p (t),K_p (t))$, to be given by
\begin{equation}\label{EqPT}
\left\lbrace\begin{array}{l}
\frac{d}{dt} \left(\begin{array}{c} V_p \\ K_p \end{array}\right)=G (V_p,K_p;V_p(t),\rho(t)) \\
V_p(0)=V_0, \quad K_p(0)=K_0,
\end{array}
\right.
\end{equation}
where $V_0$ is the initial volume of a tumor and $K_0$ its initial carrying capacity.

\subsection{Metastatic spreading}

There is no clear consensus in the biological literature about metastases being able to metastasize or not \cite{Tait2004, Bethge2012, Sugarbaker1971}. However, we argue here that cancer cells that acquired the ability to metastasize should conserve it when establishing in a new site. Moreover, since metastasis is a long process before being detectable \cite{Tait2004, Bethge2012, Sugarbaker1971} (in particular because tumors could remain dormant during large time periods), the absence of clear proof in favor of metastases from metastases could be due to the short duration of the experiments compared to the time required for a secondary generation of tumors to reach a visible size. Here we are interested in long time behaviors and, although metastases from metastases could be neglected in first approximation, we think this second order term is relevant in our setting and chose to include it in our modeling, following clinical evidences of second-generation metastases \cite{August1985}.

Successful metastatic seeding results from one malignant cell being able to overcome various adverse events including: detachment from the tumor, migration in the local micro-environment, intravasation in local blood (or lymphatic) vessels, survival in the circulation, escape from immune surveillance, extravasation, survival in a new environment and eventually establishment of a new colony at the distant site (see \cite{Gupta2006} for more details). We regroup here all these events into one emission rate $\beta(V)$ quantifying the number of successfully created metastases per unit of time, and neglect intricate description of all these processes. We assume that very small metastases do not metastasize arguing that they do not have access to the blood circulation and hence include a threshold $V_m$ below which tumors do not spread new individuals. Apart from addition of this threshold, the expression of $\beta$ is the one from \cite{iwata} and is given by 
\begin{equation}\label{EqEmission}
\beta(V)=mV^{\alpha} \mathds{1}_{V \geq V_m},
\end{equation}
where $m$ and $\alpha$ are coefficients quantifying the overall metastatic aggressiveness of the cancer disease. The parameter $m$ represents an intrinsic metastatic potential of the cancer cells. On the other hand, parameter $\alpha$ represents the micro-environment dimension of metastatic dissemination. It lies between 0 and 1 and is the third of the fractal dimension of the vasculature of the tumor under consideration. For instance, if vasculature develops superficially then $\alpha=2/3$, whereas a fully penetrating vasculature would give a value of $\alpha=1$. We chose here to take a dissemination rate only depending on the volume because simulations showed that adding a monotonous dependence on $K$ did not improve the flexibility of the model in simulations while adding at least one parameter, in opposition to our parsimony principle. 

Stating a balance law for the number of metastases when growing in size gives the first equation of \eqref{EqModel}. The boundary condition, i.e. the second equation of \eqref{EqModel}, states that the entering flux of tumors equals the newly disseminated ones. These result from two sources: spreading from the primary tumor, represented by the term $\beta(V_p(t))$ and second generation tumors coming from the metastases themselves, described by the term $\int_\Omega \beta(V)\rho(t,V,K)dVdK$. The map $\sigma \mapsto N(\sigma), \sigma=(V,K)\in \partial\Omega$ stands for the volume and carrying capacity distribution of metastases at birth.  We assume that newly created tumors have all the size $V_0$ and some initial carrying capacity $K_0$ and thus take
\begin{equation}\label{EqN}
N(\sigma)=\delta_{\sigma=(V_0,K_0)},
\end{equation}
i.e. the Dirac distribution centered in $(V_0,K_0)$ and refer to \cite{BenzekryJBD} for more detailed considerations on going from an absolutely continuous density $N$ to the Dirac mass considered here. We allow metastases to exit the domain by imposing the boundary condition only where the flux points inward. In view of expression \eqref{EqGrowthRate}, the growth velocity pushes tumors out of the domain when $K$ is less than $V_0$. This can happen when global inhibition is strong enough such that tumors can cross the line $K=V_0$, i.e. when $G_2(V_0,V_0)=bV_0-dV_0^{5/3}-eIV_0<0$. These tumors are then removed from the population, corresponding to the death of metastases caused by nutrient deprivation. 

From the solution $\rho$ of the model, relevant macroscopic quantities can be defined such as the total number of metastases: $N(t)=\int_\Omega \rho(t,V,K)dVdK$ or the total metastatic burden: $M(t)=\int_\Omega V\rho(t,V,K)dVdK$.

Overall, the model (\ref{EqModel}-\ref{EqN}) is a nonlinear transport equation endowed with a nonlocal boundary condition. The nonlinearity appears in the growth velocity $G$ through an indirect dependence on the total metastatic burden.

\subsection{Model nondimensionalization}

We are not interested here in calibrating the model parameters or outputs to relevant biological values (see \cite{Benzekry} for such a purpose) but will rather focus our interest on theoretical exploration of possible dynamics emerging from the model. Consequently, we rescale parameters and variables in order to essentialize characteristic properties of the system. Performing the following transformations - where $V^*:=\left(\frac{b}{d}\right)^{3/2}$ is the maximal reachable volume under the Hahnfeldt model \cite{dOnofrio_Gandolfi} - on variables
$$\tilde{t}=at, \; \tilde{V}=\frac{V}{V^*}, \; \tilde{K}=\frac{K}{V^*},\; \tilde{\rho}(\tilde{t},\tilde{V},\tilde{K})={V^*}^2\rho(t,V,K), \; \tilde{I}(\tilde{t})=\frac{a}{pV^*}I(t),$$
and on coefficients
$$\tilde{b}=\frac ba,\; \tilde{e}=\frac {epV^*}{a},\; \tilde{k}=\frac{k}{a}, \; \tilde{m}=\frac{m}{a}(V^*)^{\alpha}, \; \tilde{\rho^0} (\tilde{V},\tilde{K})={V^*}^2\rho^0(V,K), \; \tilde{V_0}=\frac{V_0}{V^*},\; \tilde{K_0}=\frac{K_0}{V^*},$$
gives equation \eqref{EqModel} in the new variables, with new velocity given by
\begin{equation}\label{EqGrowthRateT}
\tilde{G}(\tilde{V},\tilde {K}; \tilde{V_p},\tilde{\rho})=\left(
\begin{array}{c}
\tilde V\ln\left(\frac {\tilde{V}}{\tilde{K}}\right) \\
\tilde b\left(\tilde V- {\tilde{V}}^{2/3} \tilde K\right)-\tilde e \tilde I \tilde K
\end{array}
\right),
\end{equation}
and new differential equation on $\tilde I$
\begin{equation}\label{EqInhibT}
\frac{d\tilde I}{d\tilde t}= \tilde V_p(t)+\int_{\tilde\Omega} \tilde V\tilde\rho(\tilde t,\tilde V,\tilde K)d\tilde Vd \tilde K-\tilde k \tilde I,
\end{equation}
where $\tilde V_p =\frac{V_p}{V^*}$ and $\tilde\Omega=]\tilde V_0,+\infty[\times]0,+\infty[$. Note that this nondimensionalization made parameters $a,\,d$ and $p$ disappear from the equations. We will consider the resulting model in the following but drop the tildes to simplify the notations.

\section{Simulation results}
We present now simulations of the nondimensionalized model (\ref{EqModel}, \ref{EqPT}-\ref{EqInhibT}). For simulations based on comparison to experimental data we refer to \cite{Benzekry} and focus here on the theoretical dynamical behavior of the model. It results from the balance between three processes: dissemination (governed by the function $\beta$) growth (controlled by the growth rate $G$) and systemic inhibition of angiogenesis (represented by the inhibitor quantity $I$). For the sake of simplicity we will reduce our analysis to one parameter per process: dissemination will be represented by parameter $m$, growth by parameter $b$ and inhibition by parameter $e$. Throughout the study we fixed most parameters to $1$ (see Table \ref{TabParam}), on the notable exception of $\alpha$, fixed to $2/3$ hence assuming superficial development of the vasculature. We use this parameter set as a basis point for exploration of the parameter space. Numerical simulations were performed adapting a previously developed discretization scheme based on the straightening of the characteristics \cite{BenzekryM2AN}. 

\begin{table}[H]
\begin{tabular}{cc}
\hline
Growth & \begin{tabular}{ccc}
	  $b$ & $V_0$ & $K_0$\\
	     1 & 0.1 & 0.2\\
	 	\end{tabular} \\
\hline
SIA & \begin{tabular}{cc} 
	 $k$ & $e$ \\
	 1 & 1 \\	
		\end{tabular}\\
\hline
Dissemination & \begin{tabular}{cc}
		$m$ & $\alpha$ \\
		1 & 2/3 \\
		\end{tabular}\\
\hline	
\end{tabular}
\caption{Base parameter values. \label{TabParam}}
\end{table}

\subsection{Linear versus nonlinear dynamics}

\begin{figure}[h]
\begin{center}
\begin{tabular}{cc}
\includegraphics[width=0.4\textwidth]{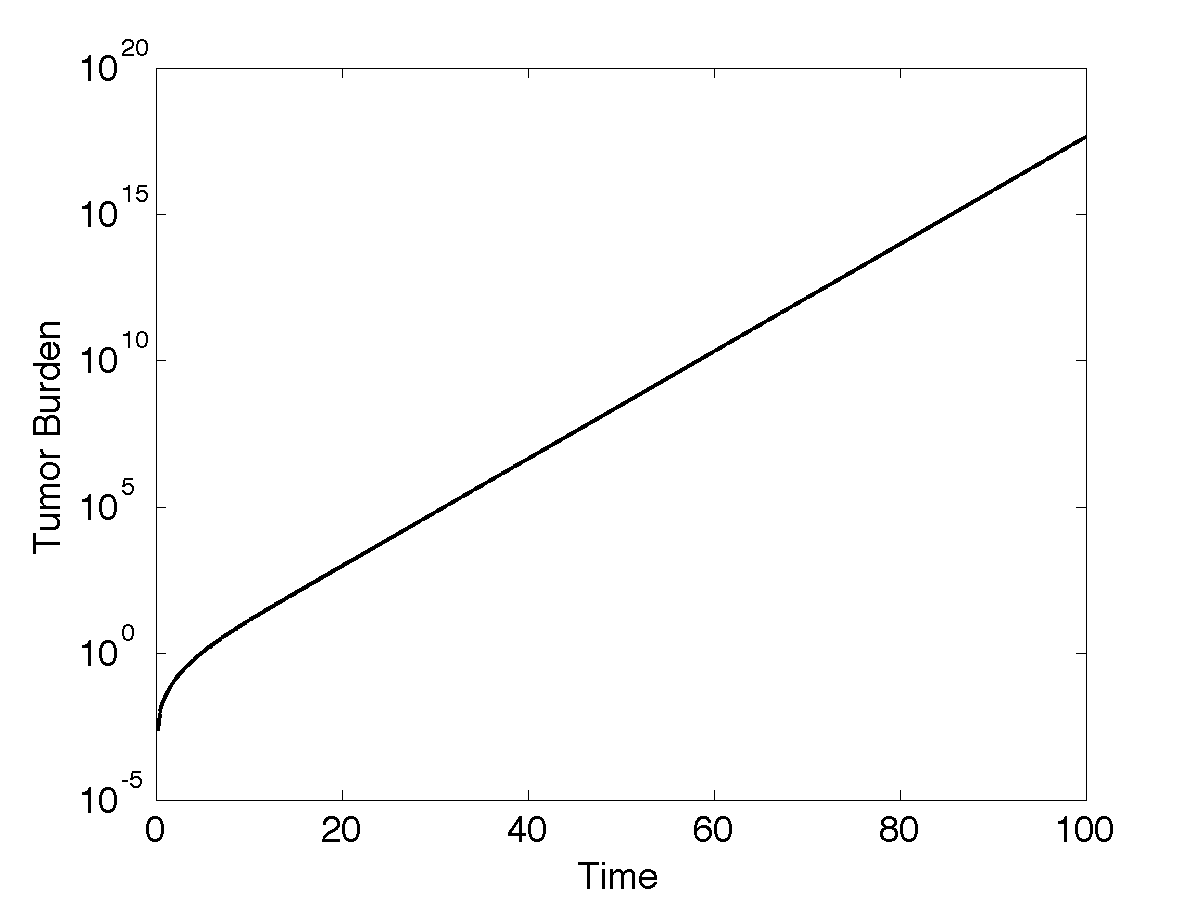} &
\includegraphics[width=0.4\textwidth]{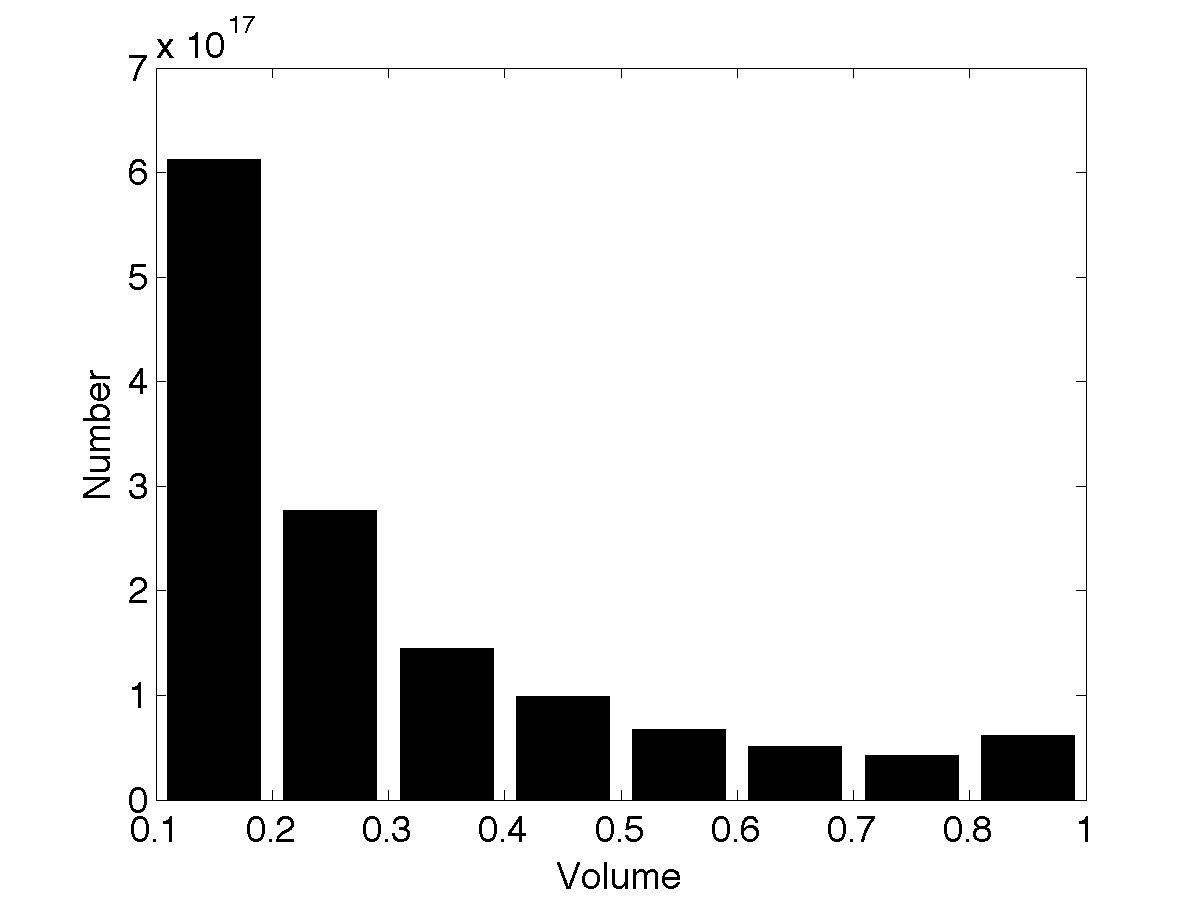} \\
Metastatic burden &  Final volume distribution
\end{tabular}
\end{center}
\caption{Asymptotic behavior of the linear model (no SIA)\label{FigLinear}. Parameters are from Table \ref{TabParam}, except $e=0$.}
\end{figure}

When SIA is not considered in the model ($e=0$), the velocity of the transport equation does not depend on $\rho$ and the model is linear. It falls then in the range of classical renewal equations from structured population dynamics (see \cite{Perthame2007} for general theory) that exhibit asynchronous exponential growth governed by the first eigenvalue  $\lambda_0$ of the underlying operator. This means that asymptotically
$$\rho(t,V,K)\sim e^{\lambda_0 t}\Pi(V,K), $$
where convergence occurs in $L^1$ norm with weight $\Psi$, the dual normalized first eigenvector, and $\Pi$ is a particular first eigenvector of the linear operator (these eigenspaces have dimension 1). The eigenvalue $\lambda_0$ is the unique solution of the following spectral equation
$$\int_0^{+\infty}\beta(X_\tau(V_0,K_0))e^{-\lambda_0 \tau}d\tau=1,$$
where $X_\tau(V_0,K_0)$ is the flow associated to the vector field $G(V,K)$. We refer to \cite{BenzekryJEE} for detailed statements and proofs of results on the asymptotic behavior of the linear model and only illustrate here in Figure \ref{FigLinear} the result and depict an example of the shape of the first eigenvector (projected on the $V$ axis).

\begin{figure}[H]
\begin{center}
\begin{tabular}{cc}
\includegraphics[width=0.4\textwidth]{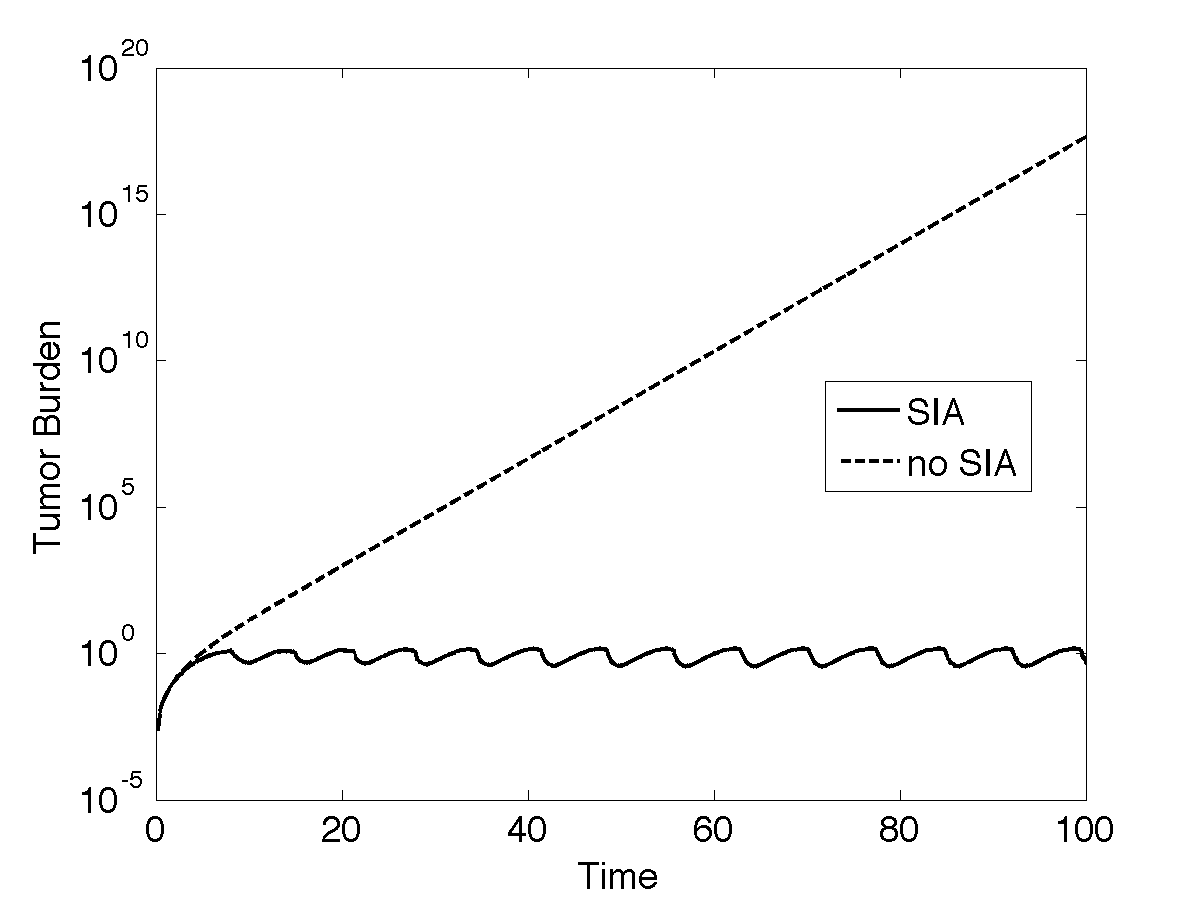} &
\includegraphics[width=0.4\textwidth]{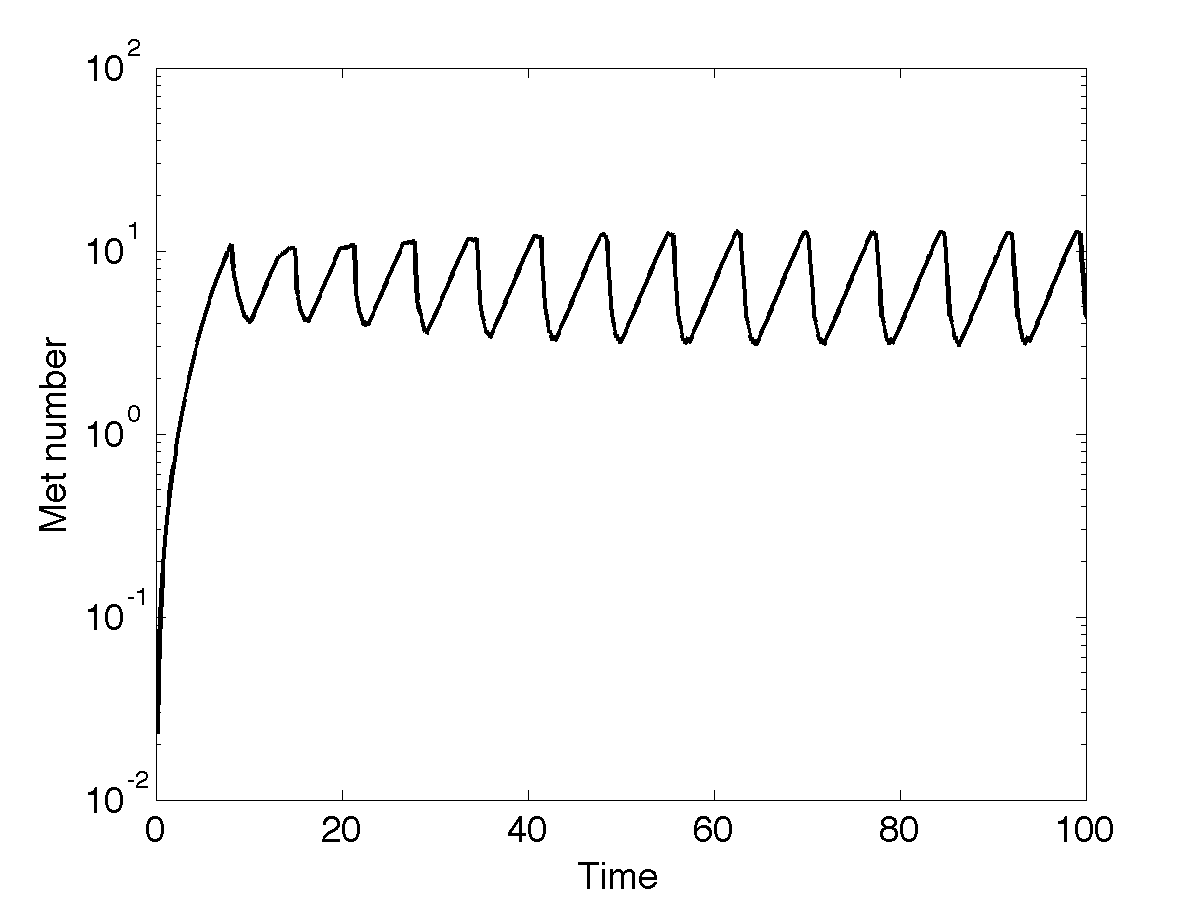} \\
Metastatic burden &
Number of metastases \\
\includegraphics[width=0.4\textwidth]{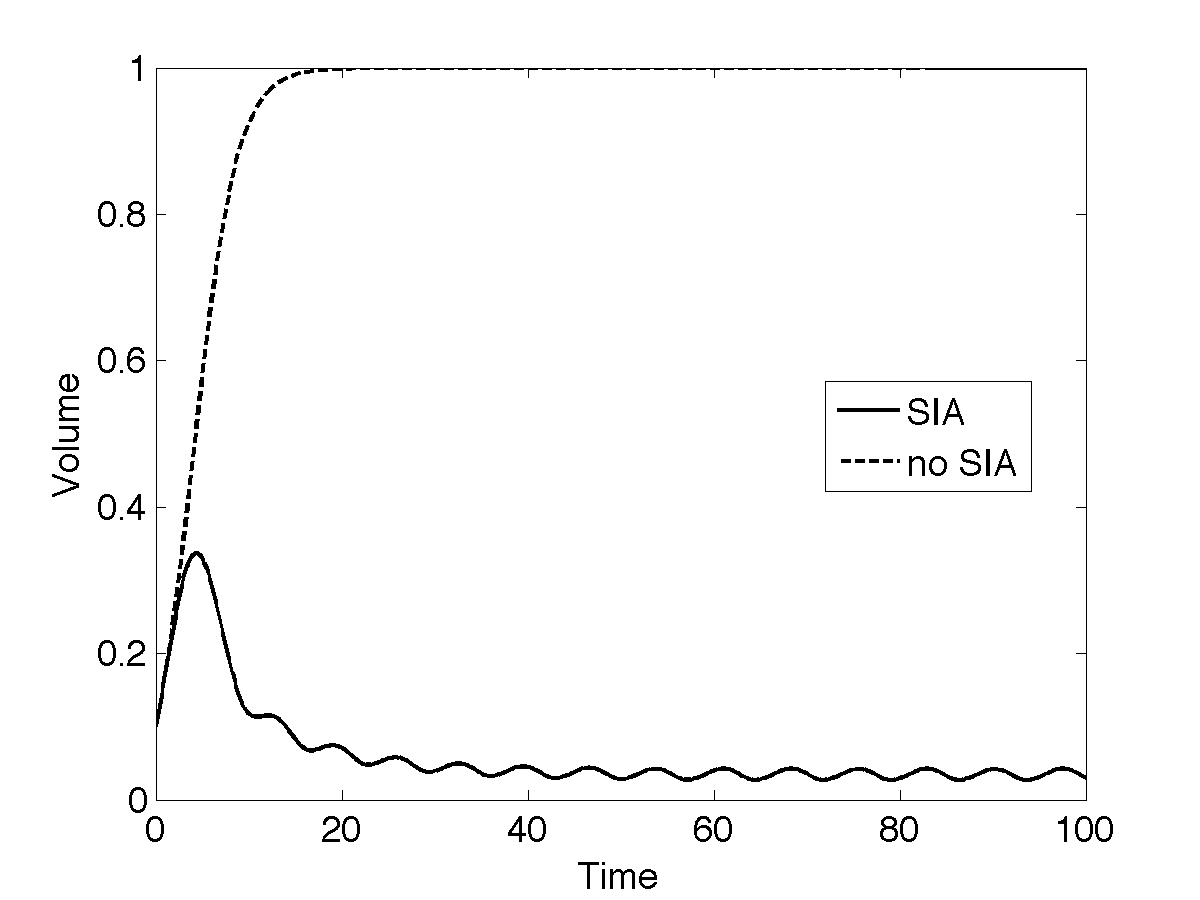} &
\includegraphics[width=0.4\textwidth]{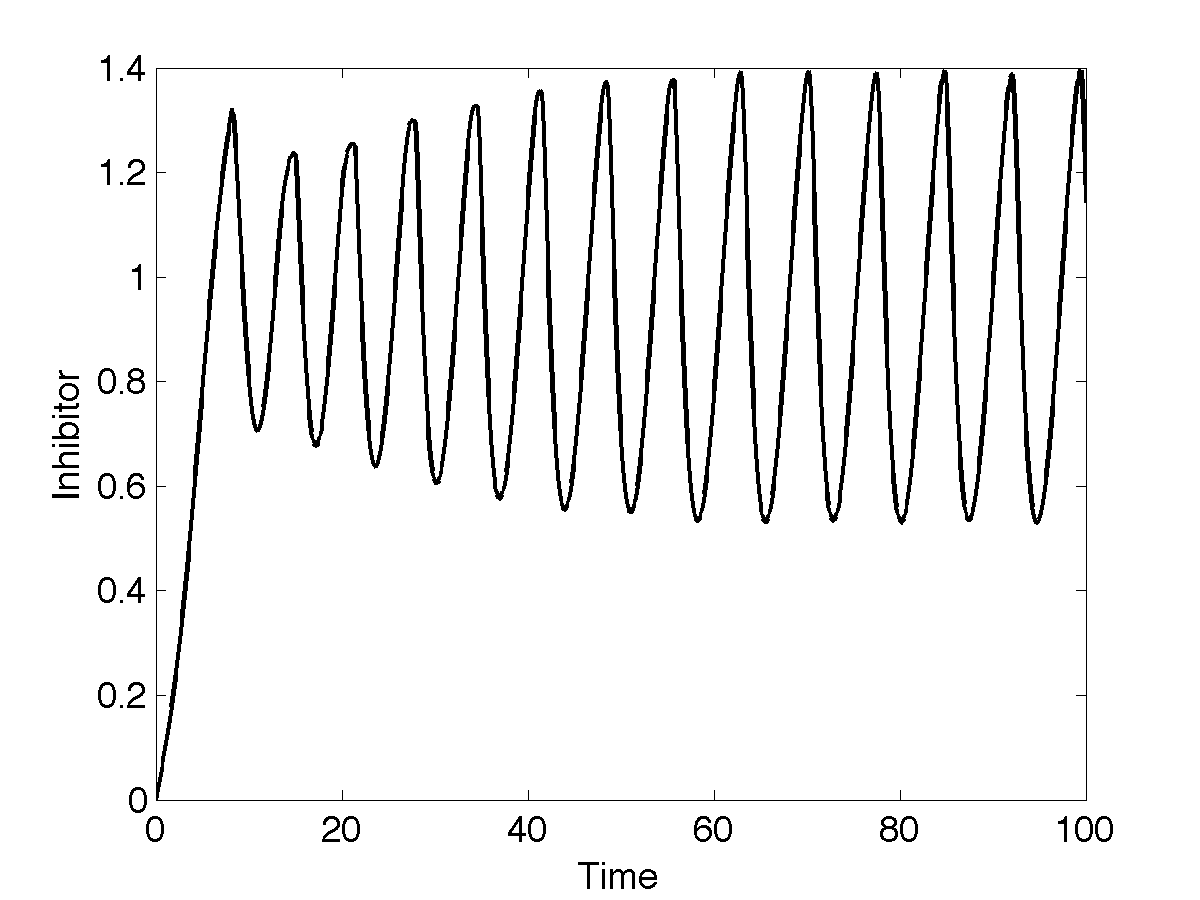} \\
Primary tumor volume &
Systemic inhibitor 
\end{tabular}
\end{center}
\caption{Dynamics of the non linear model (with SIA) \label{FigNonlinear}}
\end{figure}

When considering the nonlinear model that we introduce here, non-trivial asymptotic behavior is observed with oscillations of the system (Figure \ref{FigNonlinear}). First, under the impulsion of the primary tumor dissemination, metastases appear and start to grow until reaching a substantial metastatic burden that in turn yields growth suppression of both primary and secondary tumors. At some point, inhibition is strong enough to push metastases out of the domain, inducing a decrease in total number of metastases and metastatic burden that translates then into lower inhibitory pressure (due to clearance of the inhibitor) allowing the metastases to regrow and restart the process.

\subsection{Exploration of the growth/dissemination/inhibition balance}

We explore now the parameter space for possible different qualitative dynamics, by varying 10 fold above and below the base values of the parameters respectively controlling growth ($b$), dissemination ($m$) and inhibition ($e$). The oscillations that were observed with the base parameter, although recovered in most of the situations, are not the only possible situation as more complex dynamics are found. 

\begin{figure}[!h]
\begin{center}
\includegraphics[width=0.3\textwidth]{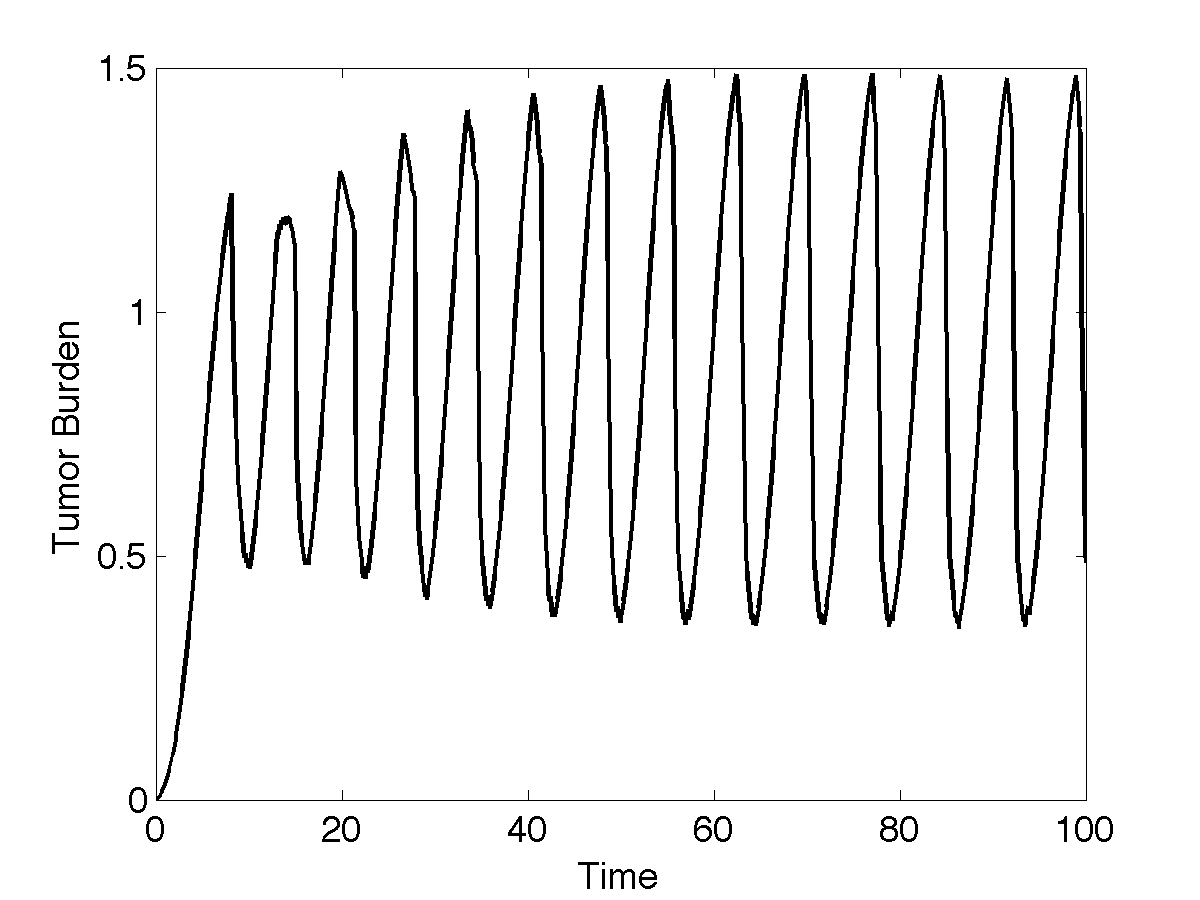}\\
Base \\
\begin{tabular}{ccc}
\includegraphics[width=0.3\textwidth]{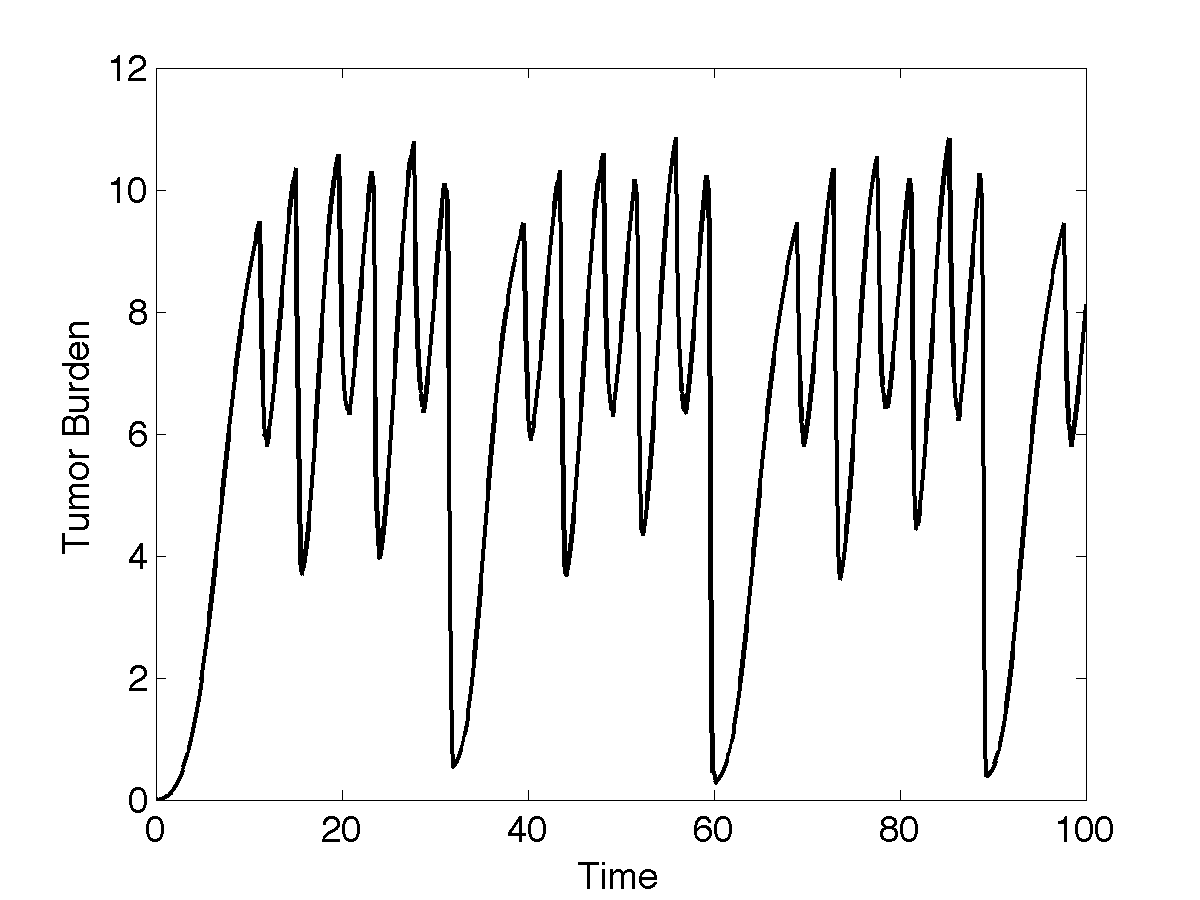} &
\includegraphics[width=0.3\textwidth]{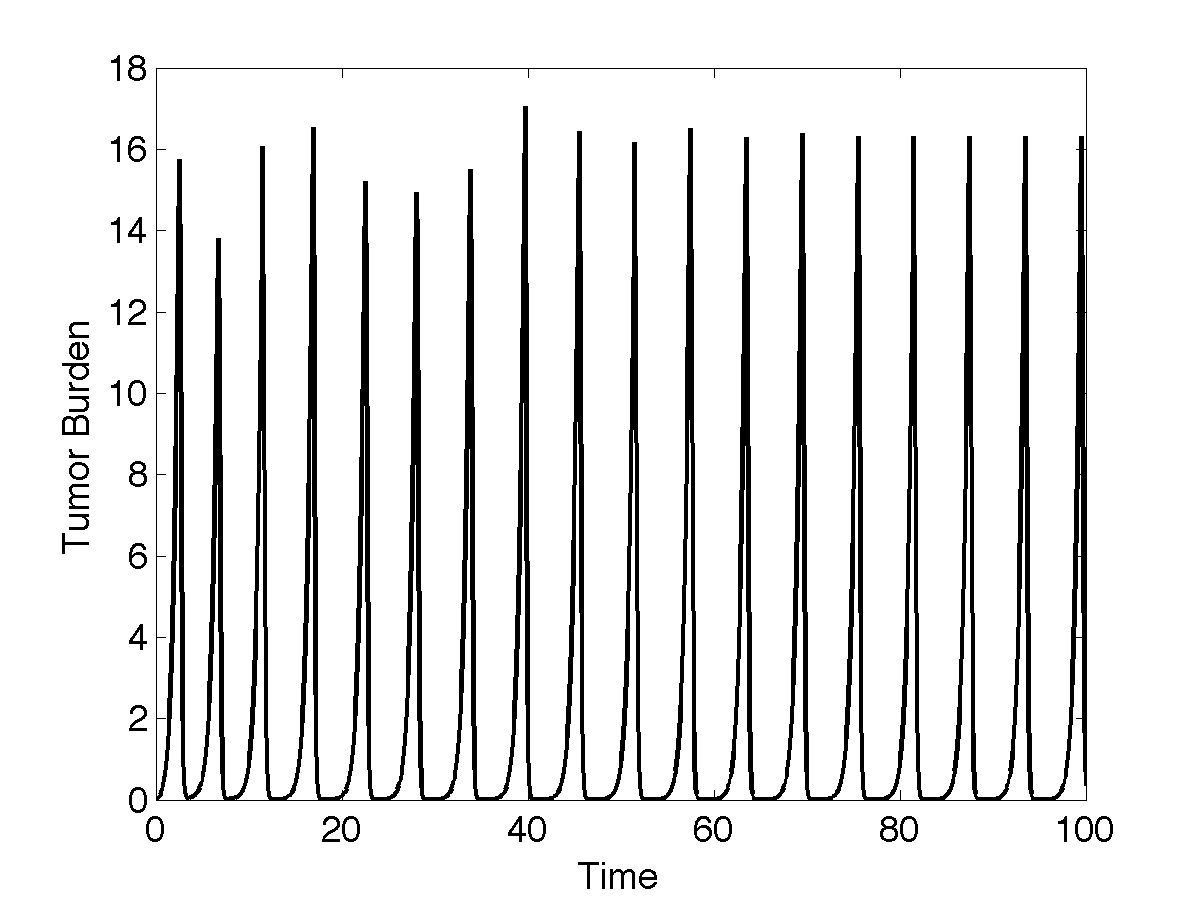} &
\includegraphics[width=0.3\textwidth]{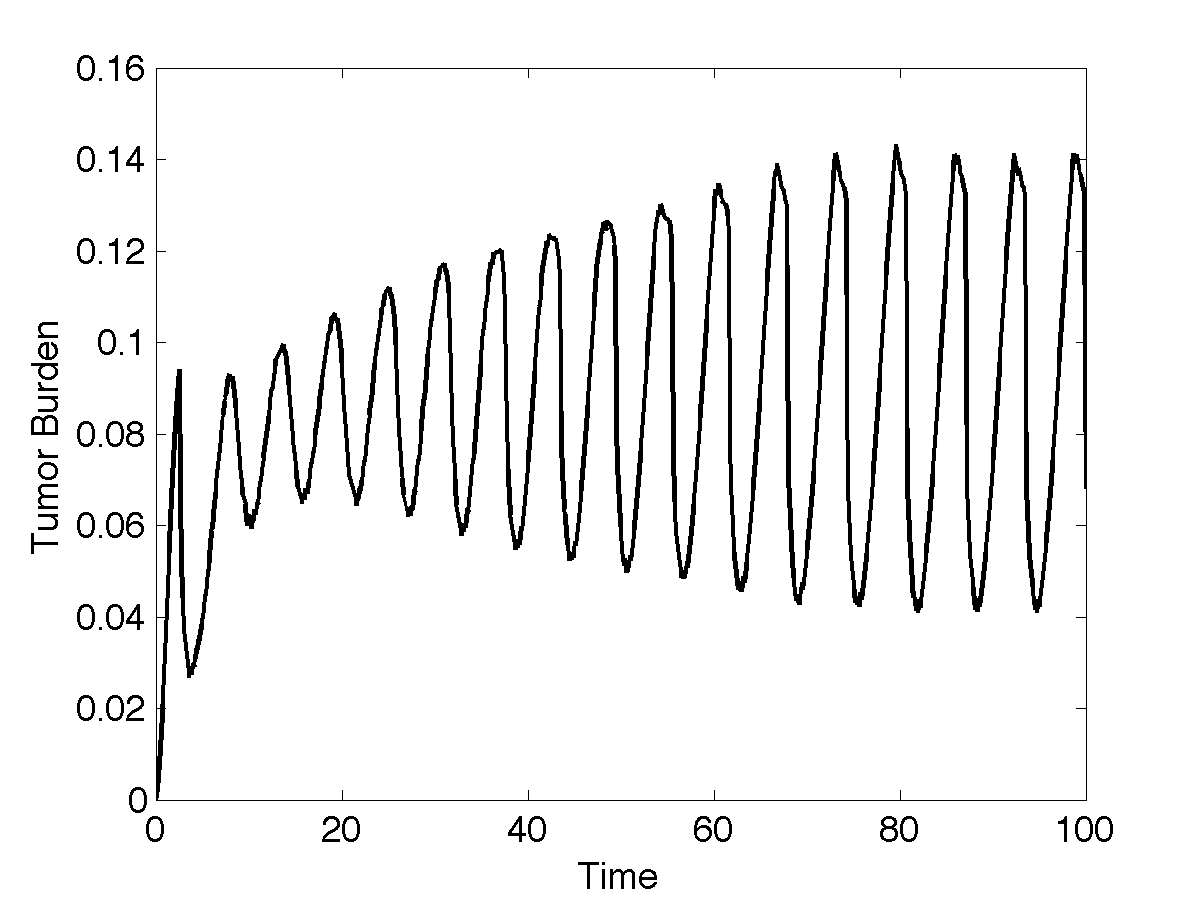} \\
Large $b$ & Large $m$ & Large $e$
\end{tabular}
\end{center}
\caption{Dynamics of the metastatic burden under 10 fold increase of representative parameters. \label{FigIncrease}}
\end{figure}

Simulation results of individual increase of each parameter are reported in Figure \ref{FigIncrease}. As appears, disruption of the base regime of parameters from Table \ref{TabParam} (where all the forces in presence are in relative equilibrium) towards more pronounced impact of either of the constitutive processes of our model generates more complex dynamics. Moreover, different parameters have different impact on the global behavior. 

Potentiation of the growth velocity (as well as resistance to the inhibition pressure) through increase of parameter $b$ results in an asymptotic behavior of the global metastatic burden which, while still being periodic, repeats a much more complex pattern, revealing interesting underlying dynamics. In particular, observation that same value of metastatic burden does not always yield same future evolution implies that no autonomous ordinary differential equation can be derived for the dynamics of $M(t)$, since same initial condition potentially leads to different future evolution. Indeed, when $M(t)$ re-reaches a previous value, the state of the global dynamical system is different because the composition of the tumors population (represented by $\rho$) is. Interestingly, this happens despite the fact that the growth rate depends on $\rho$ only through $M$ (equations (\ref{EqGrowthRateT}, \ref{EqInhibT})).

Increase of the metastatic aggressiveness of the system (parameter $m$) results in densification of the oscillations and amplitude increase, yielding sharp repeated peaks of metastatic growth. Violent increases of the total metastatic burden are followed by similarly violent decreases that make the system reach almost-zero values. 

Stronger inhibition pressure delays the stabilization of the system to an oscillatory regimen, intensifies the oscillations frequency and lowers their amplitude. Note that in this situation, as well as in the base situation, the total metastatic volume remains away from zero, suggesting a non-negligible amount of long-lasting residual disease.

\begin{figure}[H]
\begin{center}
\begin{tabular}{cc}
\includegraphics[width=0.48\textwidth]{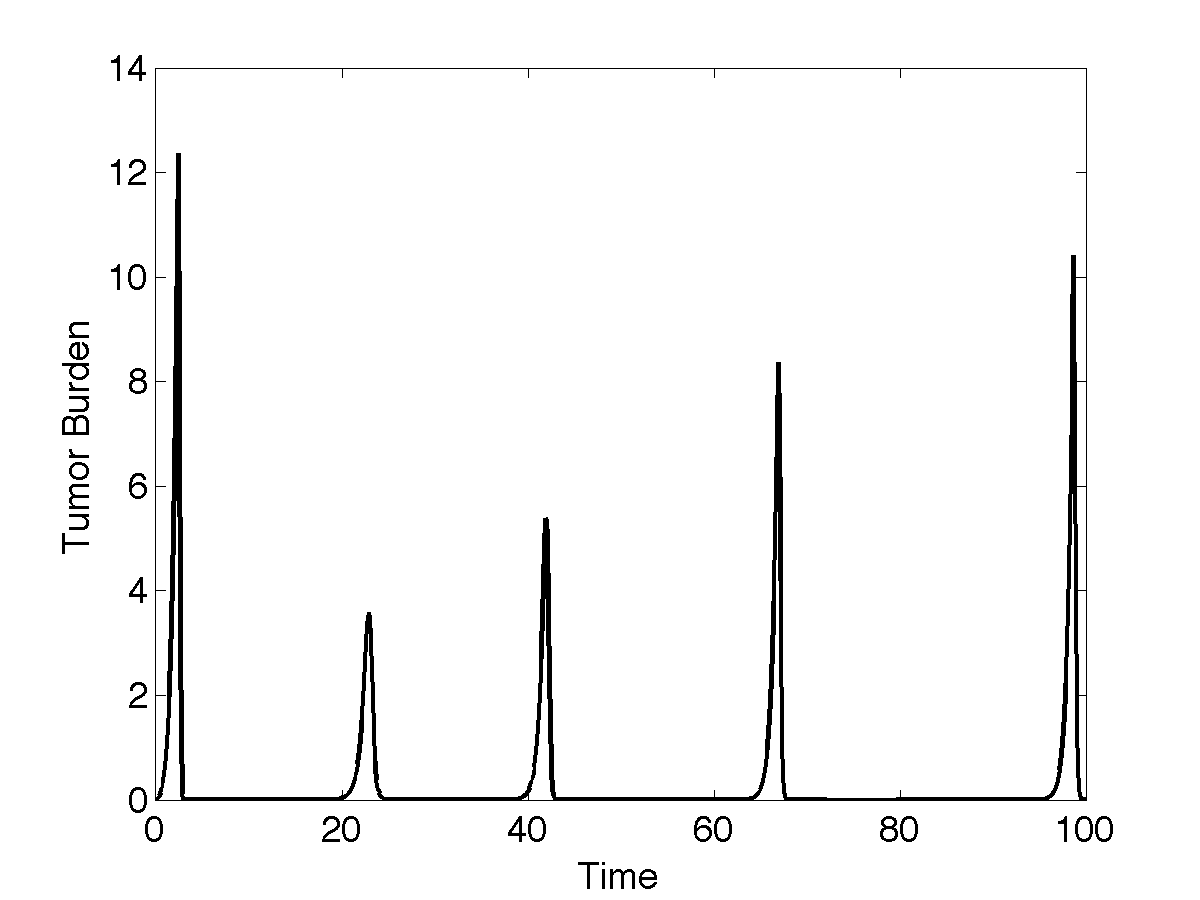}&
\includegraphics[width=0.48\textwidth]{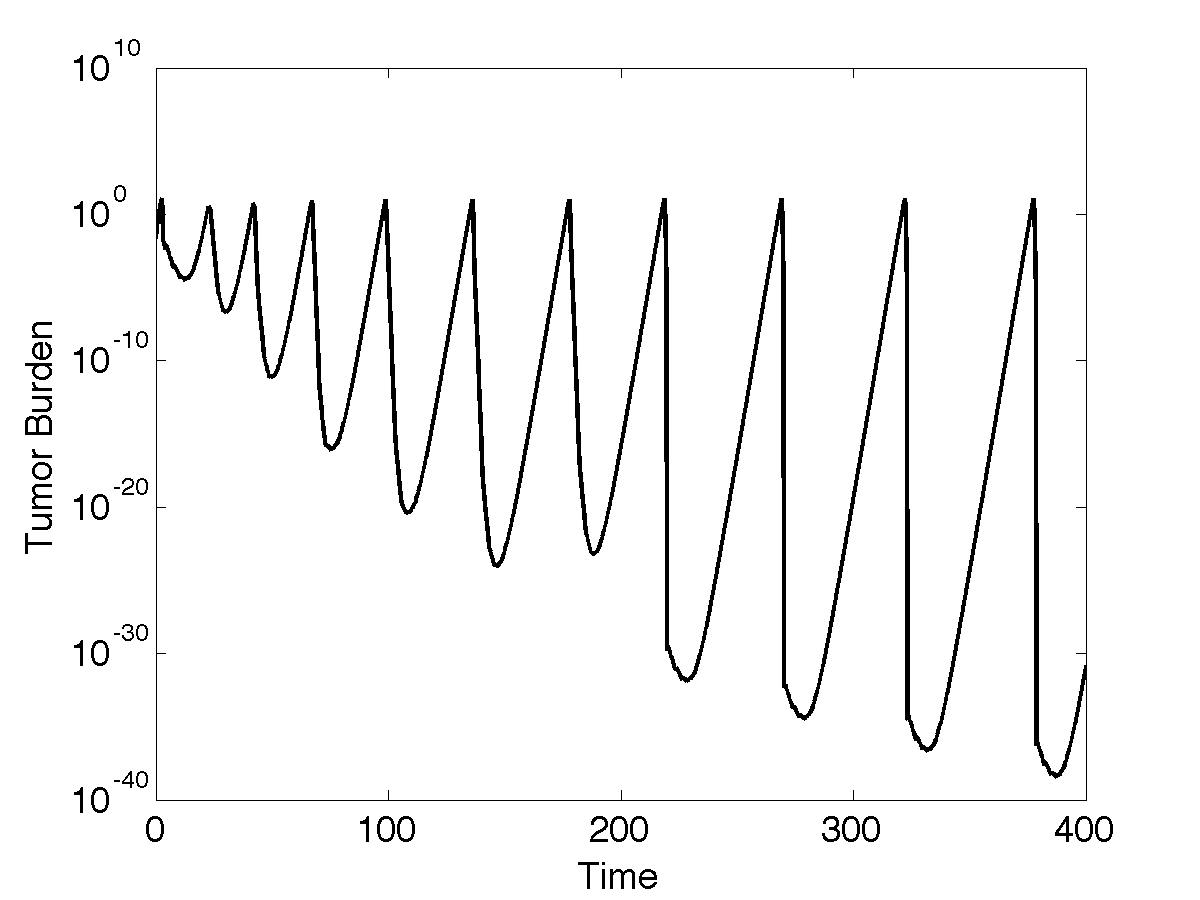}\\
$m=10, \; k=0.1$ & $m=10$, $k=0.1$  (log-scale)\\
 & \\
Increasing sharp relapses & Bounded and non-periodic \\
\includegraphics[width=0.48\textwidth]{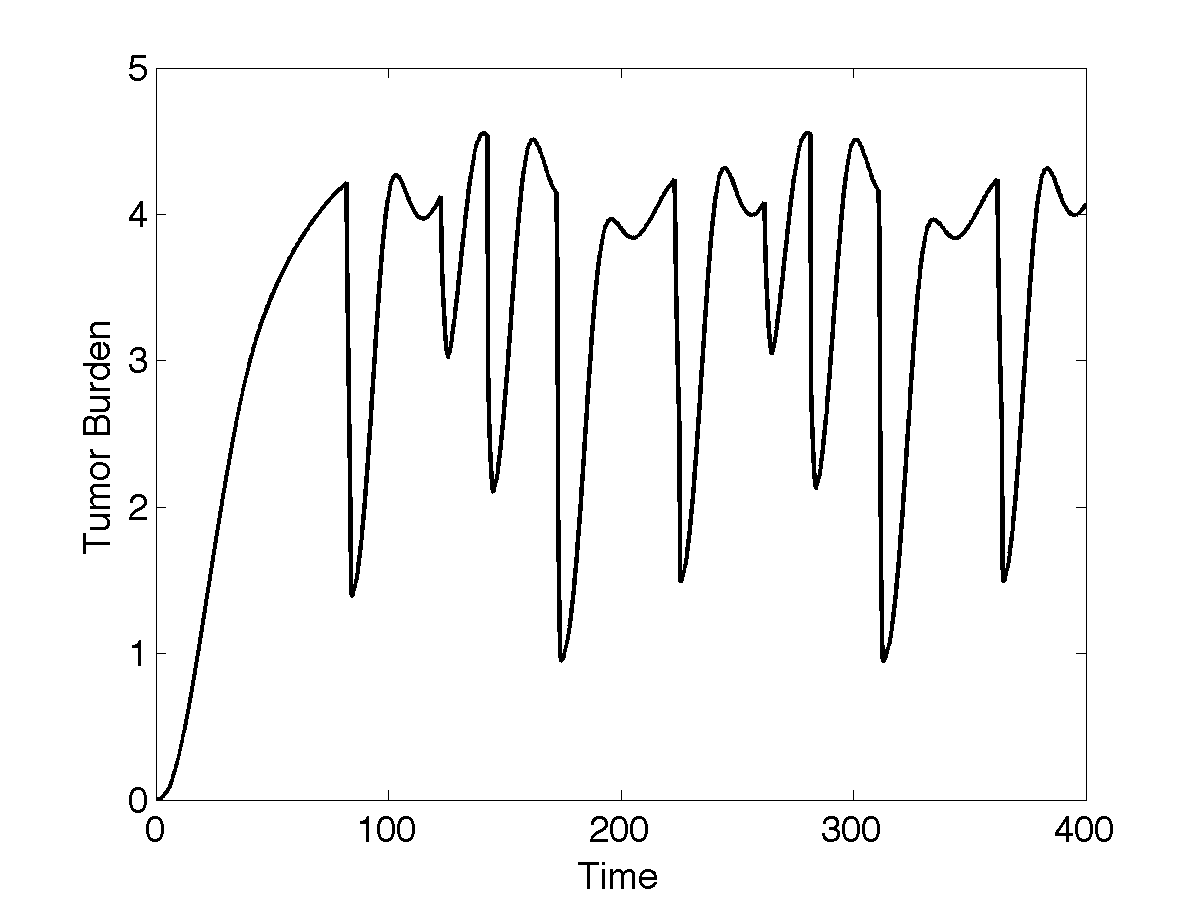}&
\includegraphics[width=0.48\textwidth]{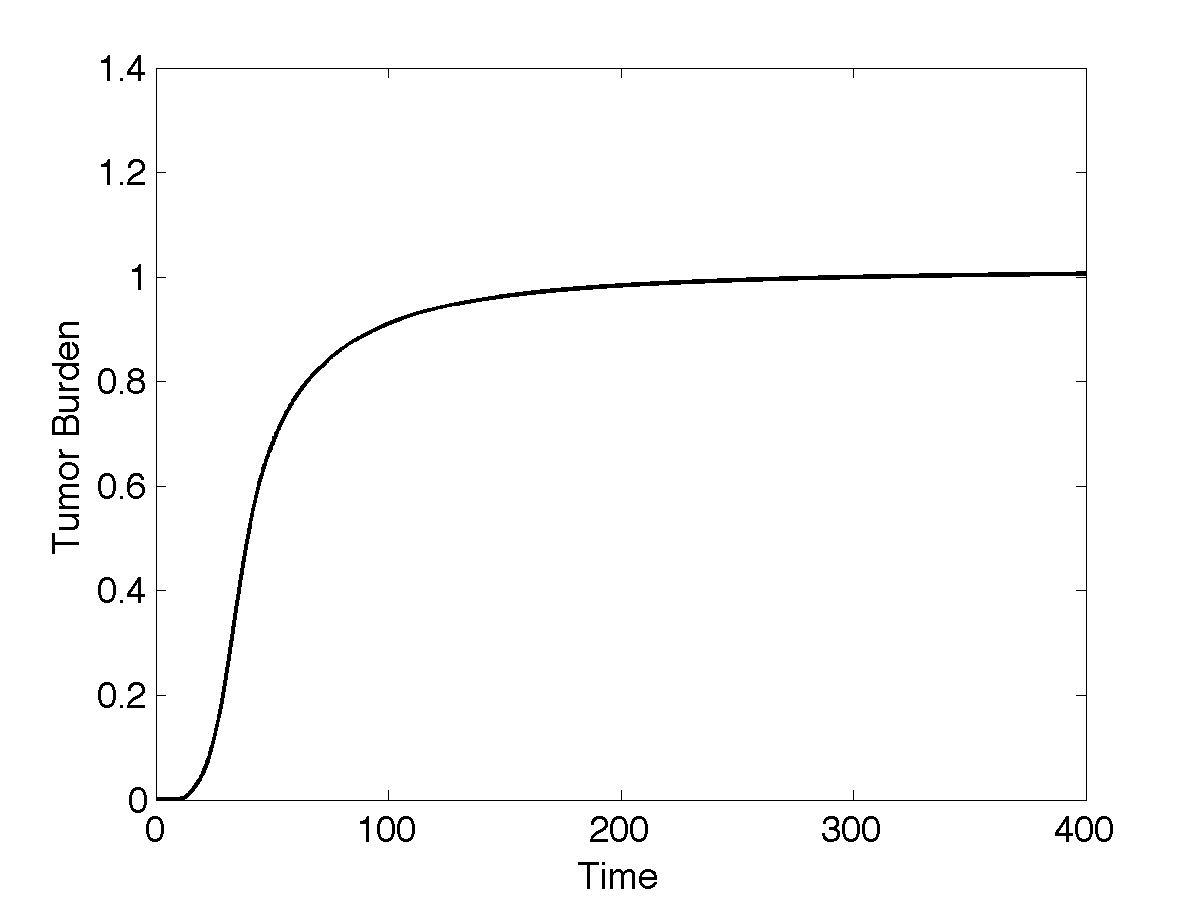}\\
 $m=0.1$, $k=0.1$, $e=0.02$ & $V_0 =10^{-4}, \quad K_0=10^{-3}$ \\
 & \\
Complex periodic behavior & Possible convergence to a stable equilibrium \\
\end{tabular}
\end{center}
\caption{Examples of other interesting dynamics\label{FigEx}}
\end{figure}

Turning our interest to the opposite situation, i.e. 10 fold decrease of the individual parameters, shows yet other interesting behaviors. Small value of $b$ generates an oscillatory asymptotic behavior with very low amplitude of the oscillations. From a biological point of view, this suggests the possibility of an homeostatic state of the system where all the forces in presence equilibrate to give a stable state where metastases don't grow while still remaining present in the organism, possibly with small volumes that would make the micro-metastases undetectable. Here, this homeostatic state stabilizes around the third of the maximal reachable size with a stable underlying distribution of metastases where volume of the largest metastasis is lower than 0.15. This result could suggest a possible explanation of reported cases of population of asymptomatic occult metastases with no evidence of primary lesion \cite{Welch2010,Folkman2004,Nielsen1987a,Black1993a} as resulting from mutual inhibitory interactions between tumors. This observation is substantiated by the simulation of the model with lower initial volume and carrying capacity (see Figure \ref{FigEx}). 

Small value of coefficient $m$ significantly delays emergence of the oscillations, since it takes more time to reach a sufficient amount of total tumor burden to trigger the dynamics. Lowering the value of $e$ has no clear impact on the frequency of oscillations but consistently increases their amplitude. Indeed, lower power of inhibitory pressure allows the metastases to grow to larger volumes.
\begin{figure}[!h]
\begin{center}
\includegraphics[width=0.3\textwidth]{BaseBurden.png}\\
Base \\
\begin{tabular}{ccc}
\includegraphics[width=0.3\textwidth]{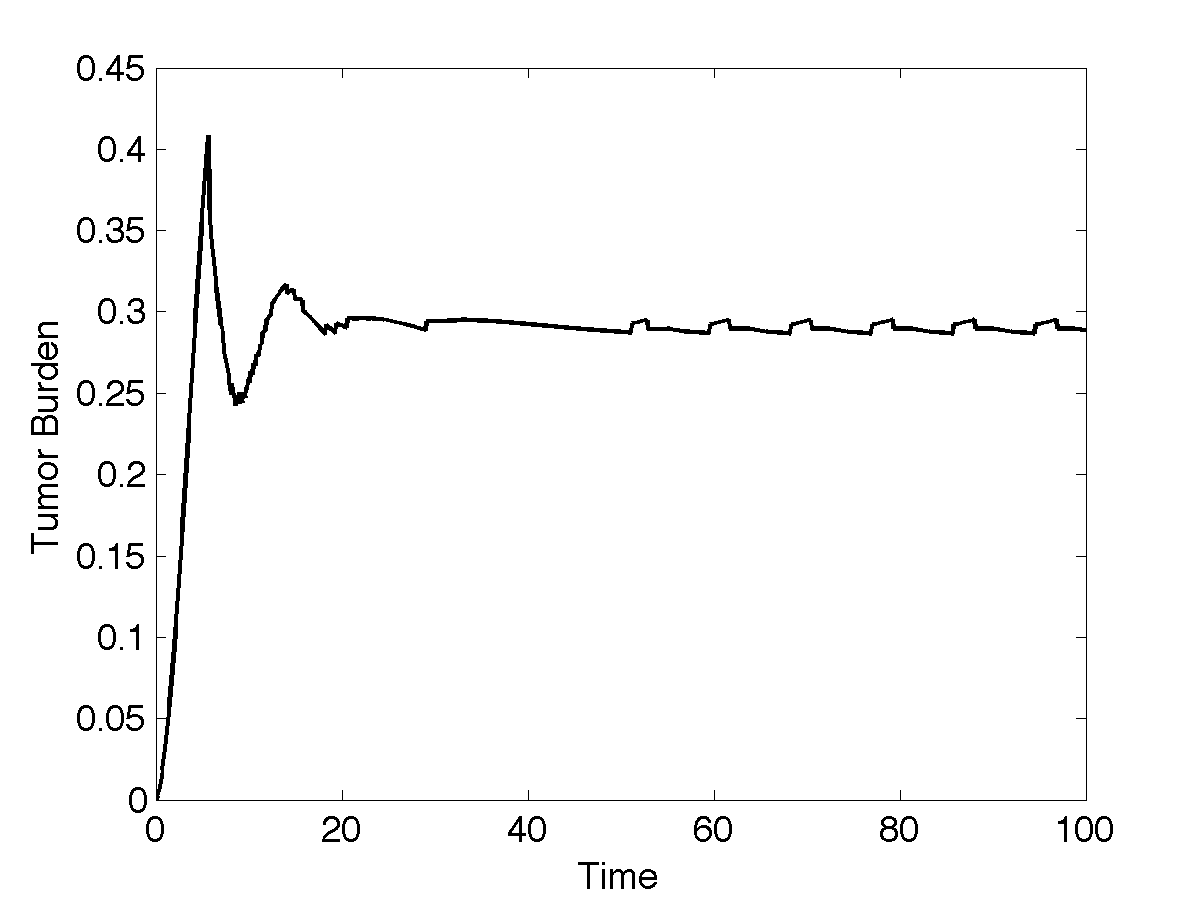} &
\includegraphics[width=0.3\textwidth]{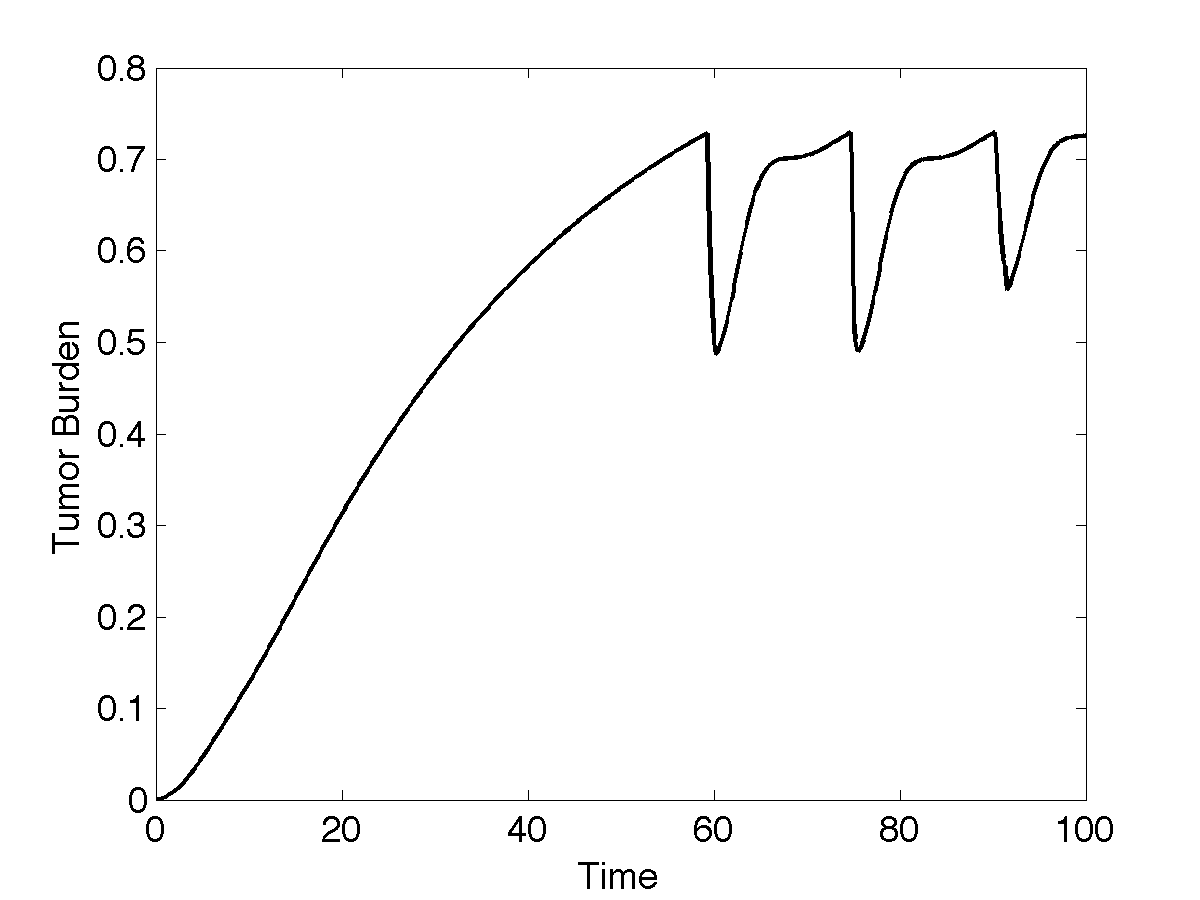} &
\includegraphics[width=0.3\textwidth]{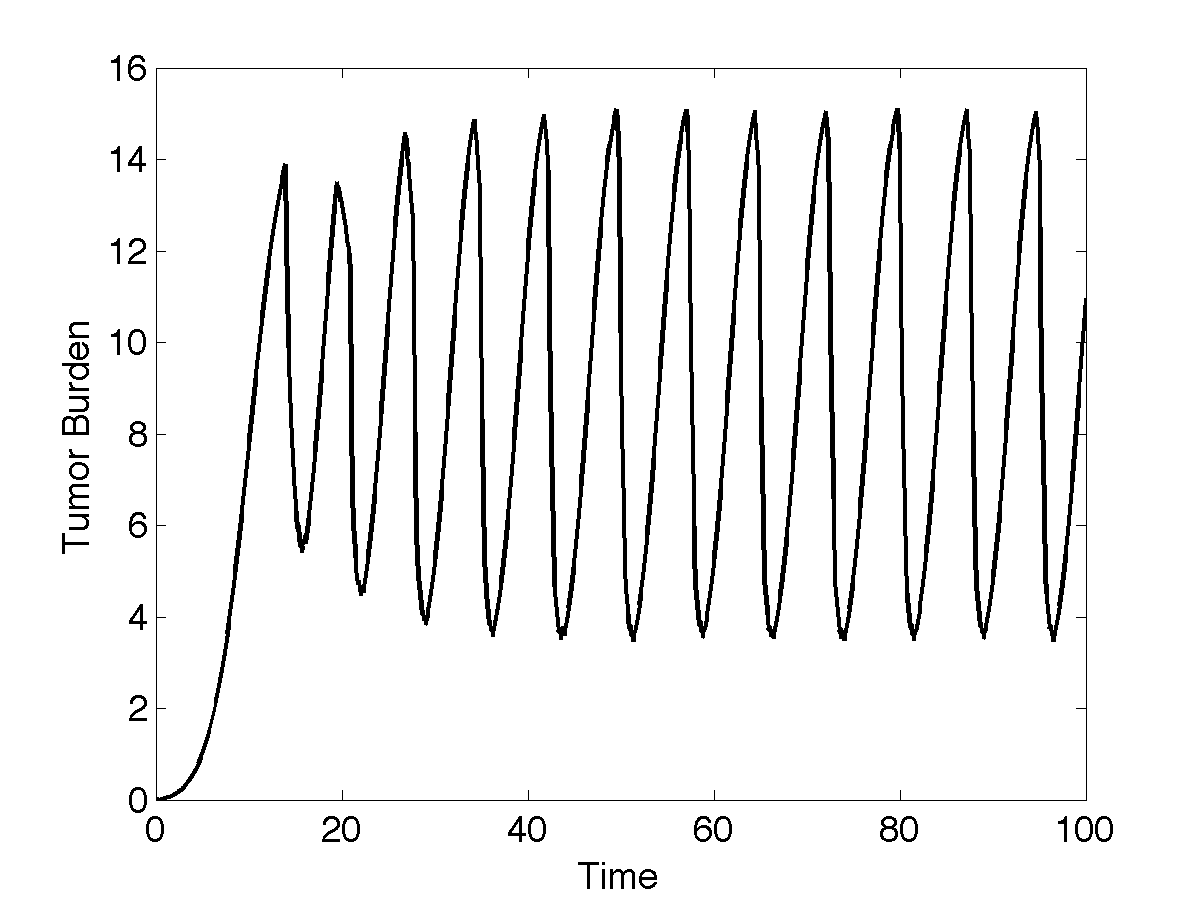} \\
Small $b$ & Small $m$ & Small $e$
\end{tabular}
\end{center}
\caption{Dynamics of the metastatic burden under 10 fold decrease of representative parameters.\label{FigDecrease}}
\end{figure}

In Figure \ref{FigEx} we report a few other examples of interesting dynamics arising from simulation of the model. In the case of large intrinsic metastatic potential ($m=10$) and inhibitor clearance set to $0.1$, after an initial sharp peak of metastatic burden, we observe firing episodes of metastatic disease of increasing intensity but delayed appearance, while rest periods are characterized by almost-zero amount of cancer mass. Although not completely relevant because of the non-biological values of the parameters, biological analogy would suggest possible violent bursts of metastatic development separated by possibly long periods of minimal (and occult) residual disease that could even lead to endogenous elimination of the cancer. Indeed, under the hypothesis of strong emission of metastases and non-negligible individual inhibitory power, it is to be expected fast exponential increase of the metastatic number of individuals, which in turn generates strong inhibition of the total population growth. However, what happens in this situation is not just acceleration of the dynamics and it remains intriguing that subsequent burst relapse with higher intensity than previous ones. A heuristic description of what happens is the following. The first burst is particular as it results from the initial condition of the primary tumor that concentrates most of the metastatic mass when it is small. This mass is eliminated all at once during the evacuation phase of the first burst that ends with smaller metastatic burden than at initiation (about $10^{-5}$ in our simulation). By the time that metastatic mass recovers significative value (say $1$ for instance), production of inhibitors has occurred that changes the condition compared to the starting point and results in lower amplitude of the second burst. At the end of this second burst, much deeper metastatic burden has been reached (of the order of $10^{-8}$) despite smaller zenith as in the first burst, because now metastases continuously outflew. With this smaller initial metastatic burden, the system had time to eliminate the inhibitor and when $M(t)$ crosses $1$ again, it does so with smaller value of $I(t)$, hence producing a higher amplitude of the relapse, which in turn provokes a smaller post-relapse burden. Repeating the same mechanism explains the following bursts. It should be noted that simulating this same situation for larger times and plotting the result in log-scale (Figure \ref{FigEx}) reveals globally bounded behavior with seemingly non-periodic orbit, thus adding another feature to the diversity of the system's dynamics.

Yet another interesting dynamics is observed for $m=0.1, \; k=0.1$ and $e=0.02$. Tormented patterns occur while still generating a periodic behavior, underlining the complexity of the dynamics of the density. On the opposite to this widely varying behavior, the model numerically exhibits convergence to a steady state for the total metastatic burden when initial conditions (for both primary tumor and metastases) are set to $(V_0, K_0)=(10^{-4}, 10^{-3})$. Same apparent convergence also occurs for the number of metastases $N(t)$ and the amount of inhibitor $I(t)$ (simulations not shown). Looking closer to the volume distribution of metastases at the end of the simulation reveals concentration of the density to the smallest possible volume, suggesting convergence to a Dirac mass located in $(V_0,K_0)$. 

\section{Conclusion}

We have presented a mathematical model for systemic inhibition of angiogenesis designed to theoretically study the relative balance of three important biological processes happening in the development of a population of metastases, namely dissemination, local growth and global (systemic) inhibition.

Simulations of this nonlinear model revealed interesting dynamics that underline the complexity. These results illustrate the interest of mathematical and computational models as useful tools for simulating the global result of a complex biology. Indeed, despite the apparent simplicity of our model that reduced here metastatic development to only three essential dynamical processes (birth, growth and inhibition), it would not be possible to guess \textit{a priori} the result of the combination of the three since, in some situations, this combination even appears as very tormented and unpredictable (see Figure \ref{FigEx}) while in others convergence to a steady state is numerically observed.

From a mathematical point of view, the numerical study we performed suggests a wide array of possible asymptotical behavior of the nonlinear model. During the exploration of the parameter space, we observed periodical behaviors with various underlying patterns, from simple oscillations to complex shapes of repeated patterns. In this context, an interesting problem would be to quantify amplitude and frequency of the oscillations in terms of the parameters of the model.

Periodicity was not always the rule since we also observed non-periodic patterns or convergence to a steady state. These observations suggest possible bifurcations in the parameter space of the infinite-dimensional dynamical system. Mathematical study of these dynamical properties, in particular how to go from a stable steady-state to a limit cycle could be interesting perspectives of our work.

From a biological viewpoint, the existence of bounded solutions of the system suggests the possibility of a stable homeostatic burden of metastases that remain in equilibrium due to mutual inhibitory interactions. In our simulations, this happens in situations where growth is substantially altered, either by a reduced growth velocity (parameter $b$) or smaller initial volume and carrying capacity of a newborn tumor with respect to what was considered here as base values (it should be noted however that within this base parameter set, $V_0$ is the tenth of the maximal reachable volume, which represents an unrealistically large initial tumor volume, making thus smaller values more biologically relevant). This modeling result could shed light on the reported observations (from forensic autopsy studies) of multiple small metastatic foci present in the organism of healthy individuals \cite{Welch2010, Nielsen1987a}, with a prevalence rate of up to 99\% (in the case of thyroid cancer) of the population, that yielded J. Folkman to use the term of ``cancer without disease" \cite{Folkman2004}. 
However, as shown in \cite{Benzekry} from calibration of the model to experimental data, elucidation of this phenomenon from systemic inhibition of angiogenesis only would require a very high production rate $p$ (or, equivalently, efficacy parameter $e$) of the systemic inhibitor. It is more likely that SIA contributes to generate such a situation in combination with other biological players (such as the immune system for instance).

\bibliographystyle{plain}

\begin{thebibliography}{}

\end{thebibliography}


\begin{thebibliography}{99}

\bibitem{Aguirre-Ghiso2007} {\sc J. Aguirre-Ghiso}. Models, mechanisms and clinical evidence for cancer dormancy. \textit{Nat Rev Cancer}, \textbf{7} (11):834–846, 2007.

\bibitem{August1985} {\sc D. A August, P. H Sugarbaker, and P. D Schneider}. Lymphatic dissemination of hepatic metastases. Implications for the follow-up and treatment of patients with colorectal cancer. \textit{Cancer}, \textbf{55}(7):1490–1494, April 1985.

\bibitem{BBHV} {\sc D. Barbolosi, A. Benabdallah, F. Hubert, and F. Verga}. Mathematical and numerical analysis for a
model of growing metastatic tumors. \textit{Math Biosci}, \textbf{218}(1):1–14, March 2009.

\bibitem{BenzekryM2AN} {\sc S. Benzekry}. Mathematical and numerical analysis of a model for anti-angiogenic therapy in metastatic cancers. \textit{ESAIM, Math. Model. Numer. Anal.}, \textbf{46}(2):207–237, 2012.

\bibitem{BenzekryJEE} {\sc S. Benzekry}. Mathematical analysis of a two-dimensional population model of metastatic growth including angiogenesis. \textit{J Evol Equ}, \textbf{11}(1):187–213, December 2011.

\bibitem{BenzekryJBD} {\sc S. Benzekry}. Passing to the limit 2D1D in a model for metastatic growth. \textit{J Biol Dynam}, \textbf{6}(sup1):19–30, January 2012.

\bibitem{BenzekryMMNP} {\sc S. Benzekry, N. Andr\'{e}, A. Benabdallah, J. Ciccolini, C. Faivre, F. Hubert and D.
Barbolosi}. Modelling the impact of anticancer agents on metastatic spreading. \textit{Math Model Nat Phenom}, \textbf{7}(1):306–336, 2012.

\bibitem{Benzekry} {\sc S. Benzekry, A. Gandolfi and P. Hahnfeldt}. Global Dormancy of Metastases due to Systemic Inhibition of
Angiogenesis. \textit{submitted}, 2013.

\bibitem{Bethge2012} {\sc A. Bethge, U. Schumacher, A. Wree and G. Wedemann}. Are metastases from metastases clinical relevant? Computer modelling of cancer spread in a case of hepatocellular carcinoma. \textit{PloS One}, \textbf{7}(4):e35689, January 2012.

\bibitem{Black1993a} {\sc W. C Black and H. G Welch}. Advances in diagnostic imaging and overestimation of disease prevalence and the benefits of therapy. \textit{N Engl J Med}, \textbf{328}(17):1237–1243, April 1993.

\bibitem{Chiarella2012} {\sc P. Chiarella, J. Bruzzo, R. P Meiss, and R. A Ruggiero}. Concomitant tumor resistance. \textit{Cancer Lett}, \textbf{324}(2):133–141, May 2012.

\bibitem{Devys2009} {\sc A. Devys, T. Goudon and P. Lafitte}. A model describing the growth and the size distribution of multiple
metastatic tumors. \textit{Discrete Contin Dyn Syst-B}, \textbf{12}(4):731–767, August 2009.

\bibitem{Dewys1972a} {\sc W. D Dewys}. Studies correlating the growth rate of a tumor and its metastases and providing evidence for tumor-related systemic growth-retarding factors. \textit{Cancer Res}, \textbf{32}(2):374–379, 1972.

\bibitem{dOnofrio_Gandolfi} {\sc A. d’Onofrio and A. Gandolfi}. Tumour eradication by antiangiogenic therapy: analysis and extensions of the model by Hahnfeldt et al. (1999). \textit{Math Biosci}, \textbf{191}(2):159–184, October 2004.

\bibitem{Fidler2003} {\sc I. J Fidler and S. Paget}. The pathogenesis of cancer metastasis: the ’seed and soil’ hypothesis revisited. \textit{Nat Rev Cancer}, \textbf{3}(6):453–458, 2003.

\bibitem{Folkman1995} {\sc J. Folkman}. Angiogenesis inhibitors generated by tumors. \textit{Mol Med}, \textbf{1}(2):120–2, January 1995.

\bibitem{Folkman2004} {\sc J. Folkman and R. Kalluri}. Cancer without disease. \textit{Nature}, \textbf{427}(6977):787, 2004.

\bibitem{Folkman1971} {\sc J. Folkman}. Tumor angiogenesis: therapeutic implications. \textit{N Engl J Med}, \textbf{285}:1182–1186, 1971.

\bibitem{Gorelik1978} {\sc E. Gorelik, S. Segal and M. Feldman}. Growth of a local tumor exerts a specific inhibitory effect on progression of lung metastases. \textit{Int J Cancer}, \textbf{21}(5):617–625, 1978.

\bibitem{Gorelik1983} {\sc E. Gorelik}. Resistance of tumor-bearing mice to a second tumor challenge. \textit{Cancer Res}, \textbf{43}(1):138–145, 1983.

\bibitem{Gupta2006} {\sc G. P Gupta and J. Massagu\ ́{e}}. Cancer metastasis: Building a framework. \textit{Cell}, \textbf{127}(4):679–695, November 2006.

\bibitem{Hahnfeldt1999} {\sc P. Hahnfeldt, D. Panigraphy, J. Folkman and L. Hlatky}. Tumor development under angiogenic signaling: a dynamical theory of tumor growth, treatment, response and postvascular dormancy. \textit{Cancer Res}, \textbf{59}(19):4770–4775, October 1999.

\bibitem{Hanahan2011} {\sc D. Hanahan and R. A Weinberg}. Hallmarks of cancer: the next generation. \textit{Cell}, \textbf{144}(5):646–674, March 2011.

\bibitem{Holmgren1995} {\sc L. Holmgren, M. S O’Reilly and J. Folkman}. Dormancy of micrometastases: balanced proliferation and apoptosis in the presence of angiogenesis suppression. \textit{Nat Med}, \textbf{1}(2):149—-153, 1995.

\bibitem{iwata} {\sc K Iwata, K Kawasaki and N. Shigesada}. A dynamical model for the growth and size distribution of multiple
metastatic tumors. \textit{J Theor Biol}, \textbf{203}(2):177–186, March 2000.

\bibitem{Nielsen1987a} {\sc M. Nielsen, J.L. L Thomsen, S. Primdahl, U. Dyreborg, and J. A Andersen}. Breast cancer and atypia among young and middle-aged women: a study of 110 medicolegal autopsies. \textit{Br J Cancer}, \textbf{56}(6):814–819, 1987.

\bibitem{O'Reilly1994} {\sc M. S O’Reilly, L. Holmgren, Y. Shing, C. Chen, R. A Rosenthal, M. Moses, W. S Lane, Y. Cao, E. H Sage and J. Folkman}. Angiostatin: a novel angiogenesis inhibitor that mediates the suppression of metastases by a Lewis lung carcinoma. \textit{Cell}, 79(2):315–28, October 1994.

\bibitem{O'Reilly1997} {\sc M. S O’Reilly, T. Boehm, Y. Shing, N. Fukai, G. Vasios, W. S Lane, E. Flynn, J. R Birkhead, B. R Olsen and J. Folkman}. Endostatin: an endogenous inhibitor of angiogenesis and tumor growth. \textit{Cell}, \textbf{88}(2):277–285, 1997.

\bibitem{Perthame2007} {\sc B. Perthame}. Transport equations in biology. \textit{Birkhauser Verlag}, 2007.

\bibitem{Prehn1993} {\sc R. T Prehn}. Two competing influences that may explain concomitant tumor resistance. \textit{Cancer Res}, \textbf{53}(14):3266–9, July 1993.

\bibitem{Retsky2010} {\sc M. Retsky, R. Demicheli, W. Hrushesky, M. Baum and I. Gukas}. Surgery triggers outgrowth of latent
distant disease in breast cancer: an inconvenient truth? \textit{Cancers}, \textbf{2}(2):305–337, March 2010.

\bibitem{Rofstad2001} {\sc E. Rofstad and B. Graff}. Thrombospondin-1-mediated metastasis suppression by the primary tumor in human melanoma xenografts. \textit{J Invest Dermatol}, \textbf{117}(5):1042–9, November 2001.

\bibitem{Sckell1998} {\sc A. Sckell, N. Safabakhsh, M. Dellian and R. K Jain}. Primary tumor size-dependent inhibition of angiogenesis at a secondary site: an intravital microscopic study in mice. \textit{Cancer Res}, \textbf{58}(24):5866–5869, 1998.

\bibitem{Sugarbaker1971} {\sc E. V Sugarbaker, A. M Cohen and A. S Ketcham}. Do metastases metastasize? \textit{Annals of surgery}, \textbf{174}(2):161–6, August 1971.

\bibitem{Tait2004} {\sc C. R Tait, D. Dodwell and K Horgan}. Do metastases metastasize? \textit{J Pathol}, \textbf{203}(1):515–8, May 2004.

\bibitem{Volpert1998} {\sc O. V Volpert, J. Lawler and N. P Bouck}. A human fibrosarcoma inhibits systemic angiogenesis and the growth of experimental metastases via thrombospondin-1. \textit{Proc Natl Acad Sci USA}, \textbf{95}(11):6343–8, May 1998.

\bibitem{Welch2010} {\sc H. G Welch and W. C Black}. Overdiagnosis in cancer. \textit{J Natl Cancer Inst}, \textbf{102}(9):605–13, May 2010.

\end{thebibliography}

\end{document}